\documentclass[aoap]{imsart}

\RequirePackage{amsthm,amsmath,amsfonts,amssymb}
\RequirePackage[authoryear]{natbib}
\RequirePackage[colorlinks,citecolor=blue,urlcolor=blue]{hyperref}
\RequirePackage{graphicx}

\startlocaldefs
\newtheorem{theorem}{Theorem}[section]
\newtheorem{lemma}[theorem]{Lemma}
\theoremstyle{remark}

\RequirePackage[colorlinks,citecolor=blue,urlcolor=blue]{hyperref}

\usepackage{graphicx}
\usepackage{tikz}
\usepackage[most]{tcolorbox}
\usepackage{pgfplots}
\pgfplotsset{compat=1.16}
\usepackage{color}
\usepackage{xcolor}
\usepackage{ulem}
\usepackage{cancel}

\endlocaldefs

\newcommand{\newtext}[1]{{\textcolor{black}{#1}}}

\newcommand{\Bin}{\mathrm{Bin}}
\newcommand{\Beta}{\mathrm{Beta}}
\newcommand{\Unif}{\mathrm{Unif}}

\newcommand{\simiid}{\overset{\mathsf{iid}}{\sim}}

\newcommand{\Ncal}{\mathcal{N}}
\newcommand{\Ucal}{\mathcal{U}}

\newcommand{\card}{\mathrm{card}}

\newcommand{\one}{\mathbf{1}}
\newcommand{\HC}{\mathrm{HC}}

\newcommand{\ex}[1]{\ensuremath{\mathbb{E}\left[ #1\right]}}

\newcommand{\goto}{\rightarrow}

\begin{document}

\begin{frontmatter}
\title{The Impossibility Region for Detecting Sparse Mixtures using the Higher Criticism}
%\title{A sample article title with some additional note\thanksref{t1}}
\runtitle{The Impossibility Region of Higher Criticism}
%\thankstext{T1}{A sample additional note to the title.}
\begin{aug}
%%%%%%%%%%%%%%%%%%%%%%%%%%%%%%%%%%%%%%%%%%%%%%
%%Only one address is permitted per author. %%
%%Only division, organization and e-mail is %%
%%included in the address.                  %%
%%Additional information can be included in %%
%%the Acknowledgments section if necessary. %%
%%%%%%%%%%%%%%%%%%%%%%%%%%%%%%%%%%%%%%%%%%%%%%
% \author[A]{\fnms{First} \snm{Author}\ead[label=e1]{first@somewhere.com}},
% \and
% \author[B]{\fnms{Third} \snm{Author}\ead[label=e2]{third@somewhere.com}}
\author{\fnms{David L.} \snm{Donoho}\ead[label=e1]{donoho@stanford.edu}},
\and
\author{\fnms{Alon} \snm{Kipnis}\ead[label=e2]{kipnisal@stanford.edu}}
%%%%%%%%%%%%%%%%%%%%%%%%%%%%%%%%%%%%%%%%%%%%%%
%% Addresses                                %%
%%%%%%%%%%%%%%%%%%%%%%%%%%%%%%%%%%%%%%%%%%%%%%
\address{
Department of Statistics,
Stanford University,
Stanford, CA 94041, USA
%\printead{e1} \printead{e2}
}
\end{aug}

\begin{abstract}
Consider a multiple hypothesis testing setting involving rare/weak effects: relatively few tests, out of possibly many, deviate from their null hypothesis behavior. Summarizing the significance of each test by a $P$-value, we construct a global test against the null using the Higher Criticism (HC) statistics of these P-values. We calibrate the rare/weak model using parameters controlling the asymptotic distribution of non-null $P$-values near zero. We derive a region in the parameter space where the HC test is asymptotically powerless. Our derivation involves very different tools than previously used to show powerlessness of HC, relying on properties of the empirical processes underlying HC. In particular, our result applies to situations where HC is not asymptotically optimal, or when the asymptotically detectable region of the parameter space is unknown. 
\end{abstract}


\end{frontmatter}



%\subtitle{{Impossibility Region} for Higher Criticism}


\maketitle

\section{Background}
\label{sec:intro} 
Multiple hypothesis testing involving rare/weak effects has attracted significant interest in statistics and machine learning over the last two decades \cite{donoho2004higher,donoho2008higher,donoho2009feature, hall2010innovated,tony2011optimal, arias2011global,cai2014optimal, arias2015sparse,li2015higher,moscovich2016exact,gontscharuk2016goodness,kipnis2019higher,DonohoKipnis2020}. Roughly speaking, the rare/weak effect setting proposes that the non-null effect is only moderately large, and, the non-null effects are present only in a relatively small, unknown subset of the tests. Here `moderately large' means that the effects are not so large that they are individually obvious; e.g. they do not individually stand out after straightforward Bonferroni correction, but they might conceivably still be detectable using more sensitive methods.

Many of the works cited above show that the Higher Criticism (HC) test of \cite{donoho2004higher} has interesting optimality properties under a phase diagram analysis in which the phase space quantifies how rare and how weak the non-null effects are. In this paper, we derive a region in the phase space where HC is asymptotically powerless, in the sense that the sum of its Type I and Type II errors approaches one within this region. 

All previous works deriving a region where HC is powerless argue that  over the given region, {\it all possible} multiple 
hypothesis tests must be asymptotically powerless. Consequently, this derivation settles matters only in cases where HC is asymptotically fully powerful inside the complement of the derived region. In such settings, HC will have the optimal phase diagram.

Actually, because the HC statistic is a practitioner's \-- rather than a theoretician's \--
tool, it can be used  in a wide variety of multiple testing situations. We have no reason to believe that it has the optimal phase diagram in general. In situations where HC does not have the optimal phase diagram, the existing strategy for
establishing its region of powerlessness does not settle matters. There remains the possibility that in certain settings, a region exists where HC is powerless, and yet some more sensitive multiple testing approach is fully powerful. 
As an example when this situation can arise, suppose that there exists a covariate predicting the presence/absence of the rare and weak effects. A naive approach is to ignore this covariate and apply HC to P-values obtained by testing each feature independently. Arguably, a better approach is to reduce each equivalence class indicated by levels of the covariate by combining the responses within this class, and applying HC to the collection of P-values obtained from each equivalent class after this reduction. As another example, suppose that many effects exist under the alternative hypothesis, but these effects are too weak for HC to be effective. Therefore, a test based on a statistics that aggregates the individual weak effects is expected to be much more powerful than HC. In Section~\ref{sec:examples} below we provide specific cases corresponding to these two examples.

\section{Problem Formulation and Statement of the Main Result
\label{sec:results}
}
 
\subsection{Rare/Weak Features}
Consider a multiple hypothesis testing  situation involving a collection of
P-values $\pi_1,\ldots,\pi_n$ associated with many independent, although not necessarily identical, tests. We are interested in testing the global null
\begin{subequations}
\label{eq:hypotheses}
\begin{align}
\label{eq:H0}
H_0~~:~~\pi_i \simiid \Unif(0,1),\quad \forall i=1,\ldots,n, 
\end{align}
against an alternative of the form 
\begin{align}
\label{eq:H1}
H_1^{(n)}~~:~~ \pi_i \sim (1-\epsilon_n)\Unif(0,1) + \epsilon_n G_{n,i},\quad \text{independently},\quad \forall i=1,\ldots,n.
\end{align}
\end{subequations}
Here $\epsilon_n = n^{-\beta}$, $\beta \in (0,1)$ is the \emph{rarity} parameter. $G_{n,i}$ is a probability distribution stochastically dominated by $\Unif(0,1)$ that specifies the behavior of the $i$th test whenever this hypothesis is non-null. 
 We assume that, for $X \sim G_{n,i}$ and $q>0$,
\begin{align}
    \label{eq:tail}
\max_{i=1,\ldots,n}\left( -\log \left(\Pr(X < n^{-q})\right) \right) \leq \log(n) \left( \alpha(q,r) + o(1) \right),
\end{align}
for some bivariate function  $\alpha(q,r)$ that is continuous, non-negative, increasing in $q$ and decreasing in $r$, and $o(1)$ indicates a sequence tending to zero uniformly in $i$ and $q\in[1/2,1]$ as $n$ goes to infinity. Note that, because $G_{n,i}$ is stochastically smaller than $\Unif(0,1)$, only the case $\alpha(q,r)\leq q$ is interesting. 
As it turns out, the large $n$ behavior of our test to be defined below only depends on the bound of the near-zero "tail" behavior of $G_{n,i}$ specified by $\alpha(q,r)$ in \eqref{eq:tail}. Effectively, \eqref{eq:tail} characterizes the way $G_{n,i}$ is smaller than $\Unif(0,1)$ on values scaling as $n^{-q}$, $q > 0$. This characterization arises in several studies of rare/weak effects considered in the past:
\begin{itemize}
    \item[(i)] The normal means model of \cite{donoho2004higher} corresponds to:
    \[
    \alpha(q,r)= (\sqrt{q}-\sqrt{r})^2.
    \]
    \item[(ii)] The two-sample normal means and the two-sample \emph{large} Poisson means models of \cite{DonohoKipnis2020} correspond to:
    \[
    \alpha(q,r)= (\sqrt{q}-\sqrt{r/2})^2.
    \]
    \item[(iii)] The two-sample \emph{small} Poisson means model of \cite{DonohoKipnis2020} corresponds to:
    \[
    \alpha(q,r) = q\frac{\log \left(\frac{2q}{ r\log(2)}\right)-1}{\log (2)} + \frac{r}{2}. 
    \]
    \item[(iv)] The normal heteroscedastic mixture model of \cite{tony2011optimal} corresponds to:
    \begin{align}
        \label{eq:alpha_mod}
    \alpha(q,r;\sigma^2)= (\sqrt{q}-\sqrt{r})^2/\sigma^2,\quad \sigma^2>0. 
    \end{align}
    More generally, \eqref{eq:alpha_mod} arises under a verity of models in which the non-null effects are on the moderate deviation scale \cite{Kipnis2021Chisquared}.
\end{itemize}

\subsection{Higher Criticism}
Let
\[
\HC_{n,i} \equiv \sqrt{n} \frac{i/n - \pi_{(i)}}{\sqrt{\pi_{(i)}\left(1-\pi_{(i)}\right)}},
\]
where $\pi_{(i)}$ is the $i$-th order statistic of $\{ \pi_i,\, i=1,\ldots,n \}$. The HC statistic is:
\begin{align}
    \label{eq:HC}
\HC_{n}^\star \equiv 
\underset{{1\leq i \leq n\gamma_0}}{\max}  \HC_{n,i},
\end{align}
where $0<\gamma_0 < 1$ is a tunable parameter that does not affect the asymptotic properties of $\HC_{n}^\star$ under $H_1^{(n)}$. It would simplify certain arguments below to fix it at
some value smaller than $1/2$, eg. $1/10$. \par
The HC test rejects $H_0$ for large values of $\HC_{n}^\star$.

\subsection{Impossibility of Detection}
We say that a sequence of tests $\{T_n\}$ based on the P-values $\pi_1,\ldots,\pi_n$ is \emph{asymptotically powerless} if 
\[
\Pr_{H_0} \left( T_{n} > h(n) \right) + \Pr_{H_1^{(n)}} \left( T_{n} \leq h(n) \right) \to 1,
\]
for any sequence $\{h(n)\}_{n\in \mathbb N}$. \par
Consider the curve $\rho(\beta)$ defined in the plane $(\beta,r)$ as 
\begin{align}
    \label{eq:rho_def}
    \rho(\beta) \equiv \sup \left\{ r>0 :\,
\max_{q\in (0,1]} \left(\frac{1+q}{2} - \alpha(q,r) - \beta\right) < 0 \right\}.
\end{align}
Our main result implies that a test based on $\HC_{n}^\star$ is asymptotically powerless when $r < \rho(\beta)$. \newtext{ Namely, whenever $(\beta,r)$ satisfies
\begin{align*}
\max_{q\in (0,1]} \left(\frac{1+q}{2} - \alpha(q,r) - \beta\right) < 0.
\end{align*}
}
\begin{theorem} \label{thm:HC_power}
Consider testing $H_0$ against $H_1^{(n)}$. Assume that $G_{n,i}$ has a continuous density $f_i$, and that, for some $C< \infty$, $\|f_i\|_{\infty} = C$ for all $i=1,\ldots,n$.
$\HC_{n}^\star$ is asymptotically powerless if $r < \rho(\beta)$.
\end{theorem}

\section{Stylized Applications
\label{sec:examples}
}
This paper's characterization of the 
region of powerlessness for HC would 
not be necessary in situations
where some lower bound argument exists
showing that no statistic -- including HC --
can be powerful over that region.
The powerlessness characterization
is much more useful in cases
where we want to show that
HC does not have the optimal
phase diagram, 
 i.e., in cases where
 there exists some {\it other} multiple testing procedure that is powerful within the region where HC is powerless. 
 In such situations, this paper's characterisation may allow
 us to establish limitations on HC which some other
procedures escape. We discuss two situations of this form.

\subsection{Coupled rare and weak effects}

Suppose that we have an array of $n \times k$ response variables $\{X_{i,j}\}_{1 \leq i \leq n, 1 \leq j \leq k}$ experiencing rare and weak effects where the noise is independent but the signal has structure: If the rare effect is present in some variable $X_{i_0,j_0}$, it is also present in the entire row $\{X_{i_0,j}\}_{j=1,\ldots,k}$. This scenario is a stylized version of the situation mentioned earlier, when a covariate predicting the presence/absence of the rare and weak effects exists. In this scenario, rather than considering each pair $(i,j)$ separately, one might instead combine data within each row to a P-value {\it per row} $i$ and then apply HC to the $n$ per-row P-values associated with each row. 
% (similar to taking norms across rows/columns and combining them: in one direction it is sparse and one direction it is dense
% worse calling out

As a concrete model, consider the following rare and weak normal means setting,
\begin{align}
\begin{split}
    H_0 \, &:\, X_{i,j}\simiid \Ncal(0,1) \quad i=1,\ldots,n,\quad j=1,\ldots,k \\
    H_1 \, & :\, (X_{i,1},\ldots, X_{i,k}) \simiid (1-\epsilon_n)\Ncal^k(0,1) + \epsilon_n \cdot 
    \Ncal^k(\mu_n,1), \qquad i=1,\ldots,n,
    \label{eq:hyp_app1}
    \end{split}
\end{align}
where $\mu_n$ and $\epsilon_n$ are calibrated to $n$ in the standard manner:
\begin{align} 
\label{eq:calibration}
    \mu_n = \sqrt{2 r\log(n)}, \qquad \epsilon_n = n^{-\beta}.
\end{align}
We compare the asymptotic power of HC in testing $H_0$ versus $H_1$ under two collections of P-values:
\begin{subequations}
    \begin{align}
        \label{eq:pval_naive}
    \pi_{i,k} & = \Pr(X_{i,j} \geq \Ncal(0,1) ),\quad i=1,\ldots,n,\quad j=1,\ldots,k,
    \end{align}
    versus
    \begin{align}
    \label{eq:pval_reduce}
    \bar{\pi}_i & = \Pr( \frac{1}{\sqrt{k}}\sum_{j=1}^k X_{i,j} \geq \Ncal(0,1) ),\quad i=1,\ldots,n. 
    \end{align}
\end{subequations}

Note that $X_{i,j} \sim (1-\epsilon_n)\Ncal(0,1) + \epsilon_n \Ncal(\mu_n,1)$ under $H_1$, hence  \eqref{eq:tail} holds for $X \sim \Ncal(\mu_n,1)$ with $\alpha(q,r) = (\sqrt{q}-\sqrt{r})^2$. The upper boundary for the of impossibility in \eqref{eq:rho_def} evaluates to
\begin{align*}
\tilde{\rho}(\beta) \equiv \begin{cases} 
\beta - 1/2 &  \frac{1}{2} < \beta < \frac{3}{4}, \\
 1-(1-\sqrt{\beta})^2 & \frac{3}{4}\leq \beta <1. 
\end{cases}
\end{align*}
It is straightforward to extend the results in Section~\ref{sec:results} to the case where the response variables are exchangeable but not necessarily independent as in \eqref{eq:hyp_app1}, implying that HC of the usual P-values \eqref{eq:pval_naive} is asymptotically powerless whenever $r < \tilde{\rho}$. On the other hand, in view of \eqref{eq:calibration} and \eqref{eq:pval_reduce}, the effect in $\bar{\pi}_i$ under $H_1$ is $k$ times stronger than the effect in $\pi_{i,j}$, while the linear reduction in the sample size from $n\cdot k$ to $n$ when moving from $\{\pi_{i,j}\}$ to $\{\bar{\pi}_i\}$ does not affect the asymptotic performance of HC. Consequently, a standard  analysis of HC applied to $\{\bar{\pi}_i\}$ implies that it has full power whenever $r > \tilde{\rho}/k$ \cite{donoho2004higher}. 

In short, the region of full power 
in the phase diagram based on reduction to
$\bar{\pi}_i$ is considerably larger than 
the region of full power for the usual 
p-values ${\pi}_i$, thereby enabling a comparison between a better phase diagram and a less good one.
% this a special case of combining across equivalent classes as mentioned before
% (similar to taking norms across rows/columns and combining them: in one direction it is sparse and one direction it is dense
% worse calling out

% Suppose that we have an array of n p-values
% where the noise is independent but the
% signal has structure, as follows.

% we can organize the p values into an n/2 by 2 array
% \pi_{i,j},  1<= i <= n/2, j=1,2

% whenever the rare weak effect is present in column j=1
% it is also present in column j=2.

% Then consider applying fisher sum_j  log(\pi_{i,j}) within each row,
% convert the fisher score to P-value, q_{i} and then
% apply HC to the new p-value q_{i} across rows. 
% let's call this HC-after-reduce

% One might guess that the phase transition boundary for
% rho(beta,HC-after-reduce) = rho(beta, HC) /sqrt(2),
% from an analogous setup in the gaussian case.

% And this seems to be so.

% This is to me a more instructive (or at least more imaginative)
% example of a situation where the rare weak model is actually
% true, but by exploiting more structure you can get an
% even  better result.

% More generally one guesses that rho(beta,HC after reduce) = rho(beta,HC)/sqrt(K) when there
% are groups of K

% More generally one guesses that if a covariate exists predicting the
%  presence/absence of rare weak effects, then applying HC-reduce across
% equivalence classes indicated by levels of the covariate, will give a better
% phase transition.

\subsection{Many weak effects but too few moderately large ones}

% Combination between dense and sparse where fisher might win
% Literature on multiple testing
Another situation where HC is ineffective compared to other tests is in the presence of a mixture of two types of effects: many weak effects and rare, moderately large effects. In certain settings it can be better to aggregate all the individual responses before reducing each response to a P-value.

As a specific model, consider the multiple hypothesis testing problem 
\begin{align}
    H_0\,:\,X_i \simiid \Ncal(0,1),\quad i=1,\ldots,n, \\
    H_1\, : X_i \simiid (1-\epsilon_n)\Ncal(a_n,1) + \epsilon_n \Ncal(\mu_n, 1),\quad i=1,\ldots,n,
\end{align}
where $\epsilon_n$ and $\mu_n$ are calibrated as in \eqref{eq:calibration}. If $a_n \to 0$ and $r < \tilde{\rho}(\beta)$, then the analysis provided in this paper implies that HC is asymptotically powerless. Nevertheless, an easy calculation shows that a test that is based on $\sum_{i=1}^nX_i^2$ is asymptotically powerful provided $a_n n^{1/2} \to \infty$.

\section{Proof of Theorem~\ref{thm:HC_power}
\label{sec:proof}
}

\subsection{Coupling}

Let $I$ be a random subset of indices in $[n] \equiv \{ 1,\ldots,n \}$ 
such that $i \in I$ with probability $\epsilon_n=n^{-\beta}$, for $i=1,\ldots,n$. %
An equivalent way of specifying $H_1^{(n)}$ of \eqref{eq:H1} sets
\begin{align}
    \pi_i \sim \begin{cases} G_{n,i} & i \in I \\
    \Unif(0,1) & i\neq I.
    \end{cases}
\end{align}
Let $Q_i^{(0)}$ and $Q_i^{(1)}$ denote the 
realized P-values under $H_0$ and $H_1^{(n)}$, respectively. 
On a common probability space, we couple these P-values so that
\begin{align}
\begin{split}
    \label{eq:Q-values}
    Q_i^{(0)} & \neq Q_i^{(1)},\quad i\in I, \\
    Q_i^{(0)} & = Q_i^{(1)},\quad i\in I^c.
\end{split}
\end{align}
For $h\in \{0,1\}$, define 
\begin{align*}
    \HC_n^{(h)} & = \max_{i=1,\ldots,\gamma_0 n} \sqrt{n} \frac{\frac{i}{n} - Q_{(i)}^{(h)} }{\sqrt{Q_{(i)}^{(h)}(1-Q_{(i)}^{(h)})}}. %= 
\end{align*}
For $h=0,1$, the joint distribution of the collection $\{Q_i^{(h)}\}_{i=1}^n$ is identical to the joint distribution of $\{\pi_i\}_{i=1}^n$ under $H_h^{(n)}$. 
{The coupling device should make clear to the reader
that the vast majority of P-values are the same actual numbers
under either hypothesis; the ones that are different appear sporadically
in $[n]$; they are concentrated in the index list $I$.
We will argue that under certain conditions, the different P-values are so thinned
out that their differences are not sufficient to allow the HC statistic to 
behave differently between the two hypotheses.}

\begin{theorem}
\label{thm:converse}
Let $\HC_n^{(0)}$ and $\HC_n^{(1)}$ be defined on a common probability space where the underlying P-values $\{Q_i^{(0)}\}_{i=1}^n$ and $\{Q_i^{(1)}\}_{i=1}^n$ are coupled as in \eqref{eq:Q-values}. Assume that $G_{n,i}$ has a continuous density $f_i$, and that, for some $C< \infty$, $\|f_i\|_{\infty} = C$ for all $i=1,\ldots,n$. Suppose that $(r,\beta)$ obeys $r < \rho(\beta)$. Then, for $c>0$,
\[
\Pr \left( \HC_n^{(1)} > \HC_n^{(0)} + c \right) \to 0,\quad \text{ as } n\to \infty.
\]
\end{theorem}
{Theorem~\ref{thm:HC_power} follows immediately.}
The proof of Theorem~\ref{thm:converse} requires 
additional notation and preliminary technical lemmas as provided below. 

% In addition, under $H_0$, it is well-known that 
% \[
%  \frac{\HC^\star_n}{\sqrt{2\log\left(\log(n)\right)}} \to_p 1.
% \]

\subsection{Preliminaries}
For any $(\beta,r) \in (1/2,1) \times (0,\infty)$, 
\begin{align}
    \label{eq:rho_cond}
\max_{q\in (0,1]} \left(\frac{1+q}{2} - \alpha(q,r) - \beta\right) < 0 \Leftrightarrow r< \rho(\beta).
\end{align}
The condition:
\begin{align}
    \label{eq:bonf_cond}
1 - \alpha(1,r) - \beta < 0 
\end{align} 
coincides with the condition $r< \rho(\beta)$ only if $q=1$ is the maximizer of \eqref{eq:rho_cond}. Therefore, any $(\beta,r)$ such that $r<\rho(\beta)$ also satisfies \eqref{eq:bonf_cond}. \par
Denote by $M \equiv M_n \sim \Bin(n,\epsilon_n)$ the RV indicating the cardinality of $I$. For $h \in \{0,1\}$, denote by
\[
F_n^{(h)}(t) \equiv \frac{1}{n} \sum_{i=1}^n \one(Q_i^{(h)} \leq t),
\] 
the empirical CDF of the P-values $\{Q^{(h)}_i\}_{i=1}^n$, and denote by
\begin{align}
    H^{(h)}(u) \equiv \frac{1}{M} \sum_{i\in I} \one(Q_i^{(h)} \leq u),
\end{align}
the empirical CDF's of P-values {in the thin subset} $I$. Define
\[
\Delta_n(t) \equiv F_n^{(1)}(t) - F_n^{(0)}(t).
\]
Also, set
\begin{align}
\label{eq:Ftilde_def}
\tilde{F}^{(h)}(t) \equiv \frac{M}{n} H^{(h)}(t).
\end{align}
Note that although $\tilde{F}^{(h)}(u)$ is scaled like the empirical CDF of $n$ random variables, there are actually only $M$ RVs involved. Furthermore, 
\begin{align}
    \Delta_n(t) = \tilde{F}^{(1)}(t) - \tilde{F}^{(0)}(t).
\end{align}
{It is crucial for our analysis that the difference CDF is normalized in some sense by $n$, 
while coupling shows that it involves only $M$ nontrivial contributions. }

{Enumerate} the elements of $I$ with indices $i_1,\ldots,i_M$ :
\[
\bar{Q}_{j}^{(h)} \equiv Q_{i_j}^{(h)},\quad i\in I. 
\]
{Our key assumption} \eqref{eq:H1} implies that
\begin{align}
    \label{eq:H1_I}
-\log(\Pr( \bar{Q}_j^{(1)} < n^{-q})) = \log(n)\left(\alpha(q,r) + o(1) \right),
\end{align}
where the term $o(1)$ is independent of $j$. Finally, consider a grid 
\begin{align}
    \label{eq:Ucal_def}
    \Ucal_{n} = \{ u_{k;n}\} \subset (0,1)
\end{align}
with spacing $n^{-3}$. 

% \subsection{Technical Lemmas}
% The proof of Theorem~\ref{thm:converse} requires the following lemmas. 

In what follows, we state and prove a series of preliminary technical lemmas. 

\begin{lemma}
\label{lem:M}
For $\beta \in(0,1)$ and $\eta>0$ such that $1-\beta+\eta>0$, define the event
\begin{align}
    \label{eq:Em_def}
    \Omega_n \equiv \{M < n^{1-\beta+\eta} \}.
\end{align}
There exists $n_0\equiv n_0(\beta,\eta)$ such that 
\begin{align}
    \Pr(\Omega_n) \leq e^{-n^{1-\beta+\eta}},\qquad n \geq n_0.
\end{align}
\end{lemma}
\subsubsection*{Proof of Lemma~\ref{lem:M}}
As $M\sim \Bin(n, n^{-\beta})$, by the Chernoff bound \cite[Ch. 4]{mitzenmacher2017probability},
\begin{align}
\label{eq:Chernof_M}
    \Pr \left( M \geq n^{1-\beta +\eta} \right) \leq \exp(- nD(n^{-\beta+\eta} \parallel n^{-\beta})),
\end{align}
where $D(x \parallel y) \equiv x \log(x/y) + (1-y) \log((1-x)/(1-y))$. 
{Inspection of the first term $x \log(x/y)$ shows
that, for some $c>0$} and $n_0 = n_0(\beta,c)$ {we have}:
\begin{align} \label{eq:DivBnd}
    n D(n^{-\beta+\eta} \parallel n^{-\beta}) \geq n^{1-\beta + \eta} \log(n) \cdot c,\quad n \geq n_0
\end{align}
for all $n \geq n_0(\beta,c)$. {Now combine  \eqref{eq:Chernof_M} with \eqref{eq:DivBnd}}. 
\qed

{The next lemma expresses a key insight.
Consider the relatively few $P$-values that are different under $H_0$ and $H_1$; 
the smallest of any of these is still much much larger 
than many $P$-values that are the same under
$H_0$ and $H_1$. Hence the ``place'' we will be looking for differences in
the behavior of $F_n^{(h)}(x)$ will not be at the far extreme $x \approx 0$.}
\begin{lemma}
\label{lem:Qminus}
Let $\bar{Q}_{-} = \min\{ \bar{Q}_{(1)}^{(0)}, \bar{Q}_{(1)}^{(1)} \}$ and $Q_{-} = \min\{ Q_{(1)}^{(0)}, Q_{(1)}^{(1)} \}$. Suppose that $(\beta,r)$ satisfies \eqref{eq:bonf_cond}. Then, for $\delta>0$, 
\begin{align}
\label{eq:lem:Qminus:1}
    \Pr \left( n \bar{Q}_- > n^{\delta} \right) \to 1,
\end{align}
and
\begin{align}
    \label{eq:lem:Qminus:2}
\Pr \left( Q_{-} < \bar{Q}_{-} \right) \to 1.
\end{align}
\end{lemma} 

\subsubsection*{Proof of Lemma~\ref{lem:Qminus}}

As $(\beta,r)$ satisfies \eqref{eq:bonf_cond} and since $\alpha(q,r)$ is continuous in a semi-neighborhood around $q<1$, there exists $\delta>0$ and $\eta>0$ such that
\begin{equation} \label{eq:assump:on:r}
 1-\alpha(1-\delta,r)- \beta-\eta < 0. 
\end{equation}
From \eqref{eq:H1_I}, we have, conditionally on $M$,
\begin{align*}
    & \Pr \left( \bar{Q}_{(1)}^{(1)} \leq n^{-1+\delta}  {|M} \right) = 
    {1 - \Pr( \bar{Q}_{j}^{(h)} > n^{-1+\delta} , j=1,\dots,M)} \\ 
    & = 1 -  \prod_{j=1}^{M} \left(1 - \Pr \left(\bar{Q}_{{j}}^{(1)} \leq n^{-1+\delta} \right) \right) \\
    & \leq 1 - \left(1 - n^{{-\alpha(1-\delta,r)}+o(1)} \right)^{M} \leq M n^{{-\alpha(1-\delta,r)}+o(1)},
\end{align*}
where in the last transition we used the inequality\footnote{{Indeed, $M \cdot \log(1+x) > \log(1+ M x)$.}} $(1+x)^M \geq 1+Mx$, $x\geq -1$; {where, of course, the term} 
 $n^{{-\alpha(q,r)}+o(1)}<1$ (ie. it represents a term which always 
belongs to $(0,1)$). 
 %On the event $\Omega_n = \{M < n^{1-\beta+\eta/2} \}$, 
 Hence unconditionally, by Lemma~\ref{lem:M}, 
 {
 \begin{align*}
 \Pr \left( \bar{Q}_{(1)}^{(1)} \leq n^{-1+\delta}  \right) 
 &\leq \Pr\left( M \cdot n^{-\alpha(1-\delta,r)+o(1)} \mid \Omega_n\right) + P( \Omega_n^c )\\
 & \leq  n^{1-\beta+\eta} \cdot n^{-\alpha(1-\delta,r)+o(1)} + e^{-n^{1-\beta+\eta/2}}.
\end{align*}
}
{Our assumption $r < \rho(\beta)$ yields \eqref{eq:assump:on:r} and, so }
 \begin{align*}
 \Pr \left( \bar{Q}_{(1)}^{(1)} \leq n^{-1+\delta}  \right) 
& \leq n^{-\eta/2 + o(1)} + o(1) \to 0;
\end{align*}
hence $\Pr \left( n \bar{Q}_{(1)}^{(1)} \geq n^{\delta} \right) \to 1$. 
Similarly, conditionally on $M$,
$\bar{Q}_{(1)}^{(0)}$ is the minimum of $M$ 
random variables $(W_j)_{j=1}^M$ (say)
iid uniform over $(0,1)$, and so: 
\begin{align*}
\Pr \left( \bar{Q}_{(1)}^{(0)} \leq  n^{\delta-1}  |M \right) & \geq 1 - \Pr \{ W_j > n^{\delta-1} ,j=1,\dots,M \} \\
& = 1 - \left(1 - n^{\delta-1} \right)^M \\
& \leq M n^{-1} \leq n^{-\beta} (1+o(1)),
\end{align*}
Unconditionally,
\[
 \Pr \left( \bar{Q}_{(1)}^{(1)} \leq n^{-1+\delta}  \right) 
 \leq n^{-\beta} (1+o(1)) + P( \Omega_n^c ) \goto 0 ,
\]
using $\beta > 0$ and (by Lemma~\ref{lem:M}) $\Pr(\Omega_n)\to 1$. 
Hence $\Pr \left( n \bar{Q}_{(1)}^{(0)} \geq n^{\delta} \right) \to 1$
and \eqref{eq:lem:Qminus:1} follows. 
For the claim in \eqref{eq:lem:Qminus:2}, 
note that $Q_{(0)}^{(0)}$ is the minimum of $n$ iid RVs uniform over $(0,1)$, hence
\begin{align*}
    \Pr\left( Q_{(0)}^{(0)} \leq \frac{\log(n)}{n} \right) = 1 -  \left(1 - \frac{\log(n)}{n} \right)^n \to 1,
\end{align*}
while \eqref{eq:lem:Qminus:1} implies
\begin{align*}
    \Pr\left( \bar{Q}_- \leq \frac{\log(n)}{n} \right) \to 0.
\end{align*}
Consequently, 
\begin{align*}
    \Pr\left( Q_{-} < \bar{Q}_- \right) & \geq \Pr\left( Q_{(1)}^{(0)} \leq  \frac{\log(n)}{n} < \bar{Q}_- \right) \\
    & \geq 1 - \Pr\left( Q_{(1)}^{(0)} > \frac{\log(n)}{n} \right) - \Pr\left( \bar{Q}_- \leq \frac{\log(n)}{n} \right) \to 1.
\end{align*}

\begin{lemma}
\label{lem:QplusB}
Under the conditions of Lemma~\ref{lem:Qminus}, we may define 
\[
Q_{+} \triangleq \max\{ Q_i^{(0)} \,:\, Q_i^{(0)} < \bar{Q}_{-} \}. 
\]
% Then 
% We have
% \[
% Q_+ = \bar{Q}_- (1 + o_p(1)),
% \]
Then for $\delta_1>0$, 
\begin{align}
\label{eq:lem:QplusB1}
n^{1-\delta_1} \left|Q_+ - \bar{Q}_- \right| \to_{p} 0.
\end{align}
Furthermore, for $\delta_2>0$, 
\begin{align}
\label{eq:lem:QplusB2}
    \Pr \left( n Q_+ > n^{\delta_2} \right) \to 1.
\end{align}
\end{lemma}

\subsubsection*{Proof of Lemma~\ref{lem:QplusB}}

Set $S \equiv \bar{Q}_{(1)}^{(0)} - Q_+$, so that $S \geq \left|\bar{Q}_- - Q_+\right|$. Now, $Q_+$ and $\bar{Q}_{(1)}^{(0)}$ correspond to a consecutive pair in  $\{Q_{(1)}^{(0)},\ldots,Q_{(n)}^{(0)}\}$. Therefore, $S \sim \Beta(1,n)$, i.e., $\Pr(S\leq t) = 1-(1-t)^n$. Consequently, for any $\epsilon>0$, 
\[
\Pr\left( n^{1-\delta_1}\left|\bar{Q}_- - Q_+\right| < \epsilon \right) \leq 
\Pr\left( S < \epsilon n^{\delta_1-1} \right) = 1 - (1 - \epsilon n^{\delta_1 - 1})^n \sim 1 - e^{-\epsilon n^{\delta_1}} \to 0.
\]
This proves the claim \eqref{eq:lem:QplusB1}. Next, fix $\delta_2>0$. From Lemma~\ref{lem:Qminus}, 
\[
\Pr\left( n \bar{Q}_- > n^{\delta_2} \right) \to 1.
\]
Take $\delta_1< \delta_2$. We have $\bar{Q}_i - Q_+ = o_p(1) n^{\delta_1-1} \le o_p(1) n^{\delta_2-1} \le o_p(1) \bar{Q}_-$, hence 
\[
Q_+ = (1+o_p(1))\bar{Q}_-,
\]
and
\begin{align*}
\Pr \left( nQ_+ > n^{\delta_2} \right) & =  \Pr \left( n\bar{Q}_- (1+o_p(1)) > n^{\delta_2} \right) \to 1,
\end{align*}
which implies \eqref{eq:lem:Qminus:2}. 
\qed

Fix $t_1<1/2$. Let $\bar{Q}_{-}\equiv \min \{\bar{Q}_{(1)}^{(0)}, \bar{Q}_{(1)}^{(1)}\}$ and $Q_{-} = \min\{Q_{(1)}^{(0)},Q_{(1)}^{(1)}\}$. Suppose that $Q_{-}<\bar{Q}_{-}$. Then the set $\{ Q_{i}^{(0)}\,:\, Q_i^{(0)} < \bar{Q}_{-} \}$ is non-empty by Lemma~\ref{lem:Qminus}. Define $Q_{+} = \max \{ Q_{i}^{(0)}\,:\, Q_i^{(0)} < \bar{Q}_{-} \}$. 

\begin{lemma}
\label{lem:HC_diff}
 We have
 {
\[
 \HC^{(1)} - \HC^{(0)}  \leq \sup_{t \in [Q_+,t_1]} \sqrt{n} \Delta_n(t) w(t),
\]
}
where 
\[
w(t) \equiv \frac{1}{\sqrt{t(1-t)}}.
\]
\end{lemma}

\subsubsection{\label{subsec:proof:lem:HC_diff} Proof of Lemma~\ref{lem:HC_diff}}

Because $Q_{(1)}^{(0)} < \bar{Q}_-$, it is also true that $Q_{(1)}^{(1)} < \bar{Q}_-$ and, moreover, that $Q_{(1)}^{(1)} = Q_{(1)}^{(0)} = Q_-$. In fact, 
\begin{align} 
\label{eq:Q_set_equality}
    \left\{ Q_i^{(0)} : Q_i^{(0)} < \bar{Q}_- \right\} = \left\{ Q_i^{(1)} : Q_i^{(1)} < \bar{Q}_- \right\} ,
\end{align}
and the {set defined this way} is nonempty. 
Also $Q^+$ is a jump point of both $F_n^{(0)}(t)$ and of $F_n^{(1)}(t)$. \par
For any {left-closed interval} $D = [a,b)$ or $D=[a,b]$, define
\begin{align}
    \HC^{(h)}_D = \max_{{t \in \{Q_i^{(h)}\} \cap D}} \sqrt{n} \cdot  {\left( F_n^{(h)}(t) - t\right)}   w(t) ; 
\end{align}
{i.e. in each case we are maximizing of the jump points of $F_n^{(h)}(t)$ in $D$}.
From \eqref{eq:Q_set_equality}, 
\begin{align*}
    \HC^{(0)}_{[Q_-,Q_+]} & =  \HC^{(1)}_{[Q_-,Q_+]} = \HC^{(0)}_{[Q_{(1)}^{(0)}, Q_+]} \\
    & = \HC^{(1)}_{[Q_{(1)}^{(1)}, Q_+]} =  \HC^{(1)}_{[Q_-, Q_+]}.
\end{align*}
Now, observe that, for $h\in\{0,1\}$, 
\[
\HC^{(h)} = \max \left( \HC^{(h)}_{[Q_{(1)}^{(h)},Q_+]}, \HC^{(h)}_{[Q_+,t_1]} \right).
\]
It follows that 
{
\[
 \HC^{(1)} - \HC^{(0)}  \leq  \HC^{(1)}_{[Q_+,t_1]} - \HC^{(0)}_{[Q_+,t_1]}.
\]
}
For a {right-open interval} $D=[a,b)$ or $(a,b)$, define
\[
\overline{\HC}_D^{(h)} = \sup_{t\in D} \sqrt{n} \left( F_n^{(h)}(t)-t \right)w(t);
\]
{i.e. the definition is maximizing over all points in $D$}.
Note that if the {left endpoint} $a$ of $D = [a,b)$ 
coincides with one of the jump points of $F_N^{(h)}(t)$, then 
\[
\HC_D^{(h)} = \overline{\HC}_D^{(h)}. 
\]
Indeed, $F_n^{(h)}(t)$ is piecewise constant with positive jumps, so that $t\to F_n^{(h)}(t) - t$ is strictly decreasing. Meanwhile, for $t \leq t_1 < 1/2$, $t \to w(t)$ is strictly decreasing. Hence, the supremum over $t\in D = [Q_i^{(h)}, t]$ always occurs at some jump point $Q_j^{(h)} \in D$. Now, for the interval $D = [Q_+,t]$, $Q_+$ is a jump point of both {$F_N^{(0)}(t)$ and $F_N^{(1)}(t)$}
{-- because of the set equality \eqref{eq:Q_set_equality}}. Hence 
\[
\HC^{(h)}_{[Q_+,t_1]} = \overline{\HC}^{(h)}_{[Q_+,t_1]},\quad h \in \{0,1\}.
\]
Finally, since $\Delta_n(t) = F_n^{(1)}(t)-F_n^{(0)}(t)$, we have
\[
 \overline{\HC}^{(1)}_{[Q_+,t_1]} - \overline{\HC}^{(0)}_{[Q_+,t_1]}  \leq \sup_{t \in [Q_+,t_1]} \sqrt{n} \Delta_n(t)  w(t).
\]
\qed

\begin{lemma}
\label{lem:spacing_U}
Let $U_{(1)}\leq \ldots \leq U_{(n)}$ denote uniform order statistics. Set $U_{(0)}= 0 $, $U_{(n+1)} = 1$, and $S_i \equiv U_{(i+1)}-U_{(i)}$, $i=0,\ldots,n$. 
{Let $\eta > 0$. 
\[
   \Pr(\min_i S_i > n^{-2-\eta}) \goto 1, \qquad n \goto \infty.
\]
also
\[
\Pr( n^{2+\eta}\min_i S_i \goto \infty ) = 1.
\]
}
\end{lemma}

\subsubsection*{Proof of Lemma~\ref{lem:spacing_U}}
We invoke familiar properties of the 
minimum of uniform spacings \cite{bairamov2010limit}, 
Set $t=x/(n(n+1))$; then for $x\geq 0$ we have
\begin{align}
& \Pr(n(n+1)\min_i S_i > x) = \Pr(\min_i S_i > t) = (1-t(n+1))^{n} = \left(1-\frac{x}{n}\right)^n  {\goto e^{-x}.} \label{eq:spacing_proof}
\end{align}
{Both conclusions follow directly.}
%\nr{Transition above is unclear. If $S_i \sim \Beta(1,n)$, shouldn't we replace $(1-(n+1)t)^n$ with $(1-t)^{n(n+1)}$}
\qed

\begin{lemma}
\label{lem:E1}
Consider the grid $\Ucal_{n}$ of \eqref{eq:Ucal_def}. For $U_1,\ldots,U_n \simiid \Unif(0,1)$, let $\Omega_n^{(1)}$ denote the event that there exists, for $i=1,\ldots,n$, at least one grid point situated between each consecutive pair $[U_{(i)}, U_{(i+1)}]$, $1\leq i \leq n-1$. Then 
\[
\Pr(\Omega_n^{(1)}) \to 1.
\]
\end{lemma}

\subsubsection*{Proof of Lemma~\ref{lem:E1}}
Lemma~\ref{lem:spacing_U} implies that the maximal spacing between each consecutive pair 
eventually exceeds $n^{-2-\eta}$, for each $\eta>0$. 
Because the spacings in $\Ucal_n$ all equal $n^{-3}$, 
the probability that there exists at least one grid point between each consecutive 
pair tends to one.
\qed

\begin{lemma}
\label{lem:E2a}
Consider $X_i \sim P_{X_i}$, for $i=1,\ldots,n$ independently, where we assume that each $P_{X_i}$ has a continuous density function $f_{X_i}$, and that $\|f_{X_i}\|_{\infty} = C < \infty$ for all $i=1,\ldots,n$ for some $C>0$. Consider the grid $\Ucal_{n}=\{u_k\}$ of \eqref{eq:Ucal_def}. Set $T_i \equiv X_{(i+1)}-X_{(i)}$. {Then for $\eta > 0$,
\begin{align}
    \Pr \left\{ n^{2+\eta} \min_i T_i > x \right\} \goto 1, 
\end{align}
and
\[
\Pr( n^{2+\eta}\min_i T_i \goto \infty ) = 1.
\]
}
\end{lemma}

\subsubsection*{Proof of Lemma~\ref{lem:E2a}}
{Let $F_X = \mathrm{Ave}\{F_{X_i}\}_{i=1}^n$ denote the marginal distribution;
it has {continuous} density $f_{X}$ obeying $ f_X \leq C$ as well.}
Now define RV's $U_{{(i)}}$ by
{
\[
F_{X}(X_{i}) = U_i, \quad i=1,\ldots,n.
\]
These are iid uniform RV's obeying 
\[
{F_{X}(X_{(i)}) = U_{(i)}, \quad i=1,\ldots,n.}
\]
}
Set $X_{(0)} = 0$, $X_{(n+1)}=1$
The
RV's $T_i = X_{(i)} - X_{(i-1)}$, {for $i=1,\dots,n+1$ may be written as}
\[
T_i = \frac{U_{{(i)}}-U_{{(i-1)}}}{f_X(x_i^*)}, \qquad =1,\dots,n+1,
\]
where $x_i^*$ is chosen {by the mean value theorem} to obey
\[
U_{{(i+1)}} - U_{{(i)}}  = F_X(X_{(i+1)}) - F_X(X_{(i)}) = f_X(x_i^*) \left(X_{(i+1)}-X_{(i)} \right).
\]
Because $f_X(x_i^*) \leq C$, 
\[
T_i = X_{(i)} - X_{(i-1)} \geq \frac{1}{C} \left( U_{{(i+1)}}-U_{{(i)}} \right). 
\]
Set {$S_i = U_{{(i+1)} }-U_{{(i)}}$}, $i=1,\dots, n$.
%and $S_i = U_{(i+1)}-U_{(i)}$, where $U_{(i)}$ is according to the true ordering of $U_1,\ldots,U_n$. 
Lemma~\ref{lem:spacing_U} implies {
\begin{align}
    \Pr \left\{ n^{2+\eta} \min T_i \geq x \right\} & \geq 
     \Pr \left\{ C n^{2+\eta} \min S_i \geq x \right\} \goto 1. 
\end{align}
}

\begin{lemma}
\label{lem:E2}
Let $\Omega_n^{(2)}$ denote the event that there exists, for $i=1,\ldots,n$, at least one grid point situated between each consecutive pair $[Q^{(1)}_{(i)},Q^{(1)}_{(i+1)}]$, $1\leq i \leq n-1$. Then $\Pr(\Omega_n^{(2)})\to 1$
\end{lemma}

\subsubsection*{Proof of Lemma~\ref{lem:E2}}
Using Lemma~\ref{lem:E2a} with 
$X_i = Q_i^{(1)}$, for $i=1,\ldots,n$, and 
$x=\log(n)$, there exists $n_0$ such that 
\begin{align}
    \Pr \left\{ C n(n+1) \min T_i \geq \log(n) \right\} < 1/n,\quad n \geq n_0.
\end{align}
Consequently, the probability that the minimal spacing between each consecutive pair among $\{Q_i^{(1)},\,i=1,\ldots,n\}$ is greater than $n^{-2.5}$ is tending to 1. The stated claim follows because the spacings in the grid $\Ucal_n$ all equal $n^{-3}$.

\begin{lemma}
\label{lem:2a}
Let $F(x)$ denote the empirical CDF of $x_1,\ldots,x_n \subset (0,1)$, 
and let $\Ucal_{n}=\{u_k\}$ be the grid of \eqref{eq:Ucal_def}. Assume that there exists, for $i=1,\ldots,n-1$, at least one grid point situated between each consecutive pair $(x_{(i)},x_{(i+1)})$. Then for all sufficiently large $n$,
\begin{align*}
    \max_{i=1,\ldots,n} \sqrt{n} \frac{F(x_i)-x_i}{\sqrt{x_i(1-x_i)}} - \max_{u_k \in \Ucal_n} \sqrt{n} \frac{F(u_{k})-u_{k}}{\sqrt{u_{k}(1-u_{k})}} \leq n^{-1}.
\end{align*}
\end{lemma}

\subsubsection*{Proof of Lemma~\ref{lem:2a}}

\begin{proof}
{Let $u_{k(i)}$ denote the smallest value of $u_k$ satisfying $ u \in (x_{(i)},x_{(i+1)})$
and $u \in \Ucal_{n}$.  
Then 
\begin{equation} \label{eq:matchCDF}
    F_n^{(1)}(x_i)=F_n^{(1)}(u_{k(i)}), \qquad i=1,\dots,n .
\end{equation}
Hence:
\begin{align*}
\max_{i=1,\ldots,n} & \left( \frac{F(x_i)}{\sqrt{x_i(1-x_i)}} -  \frac{F(u_{k(i)})}{\sqrt{u_{k(i)}(1-u_{k(i)})}} \right) \\
= &
\max_{i=1,\ldots,n} F(x_i) \left( \frac{1}{\sqrt{x_i(1-x_i)}} -  \frac{1}{\sqrt{u_{k(i)}(1-u_{k(i)})}} \right)  \\
\leq&
\max_{i=1,\ldots,n} \left| \frac{1}{\sqrt{x_i(1-x_i)}} -  \frac{1}{\sqrt{u_{k(i)}(1-u_{k(i)})}} \right|. 
\end{align*}
Also
\[
\max_{u_k \in \Ucal_n}  \frac{F(u_{k})-u_{k}}{\sqrt{u_{k}(1-u_{k})}} \geq 
\max_{i=1,\ldots,n}  \frac{F(u_{k(i)})-u_{k(i)}}{\sqrt{u_{k(i)}(1-u_{k(i)})}},
\]
admittedly, generally with equality. Combining these:
\begin{align*}
\max_{i=1,\ldots,n} & \frac{F(x_i)-x_i}{\sqrt{x_i(1-x_i)}} - \max_{u_k \in \Ucal_n}  \frac{F(u_{k})-u_{k}}{\sqrt{u_{k}(1-u_{k})}} \\
\leq &\max_{i=1,\ldots,n} \left( \frac{F(x_i)-x_i}{\sqrt{x_i(1-x_i)}} -  \frac{F(u_{k(i)})-u_{k(i)}}{\sqrt{u_{k(i)}(1-u_{k(i)})}} \right) \\
\leq & \max_{i=1,\ldots,n}  \left| \frac{F(x_i)}{\sqrt{x_i(1-x_i)}} -  \frac{F(u_{k(i)})}{\sqrt{u_{k(i)}(1-u_{k(i)})}} \right| \\
+ & \max_{i=1,\ldots,n} \left|\frac{x_i}{\sqrt{x_i(1-x_i)}} - \frac{u_{k(i)}}{\sqrt{u_{k(i)}(1-u_{k(i)})}}\right| \\
\leq &  \max_{i=1,\ldots,n} \left|\frac{1}{\sqrt{x_i(1-x_i)}} - \frac{1}{\sqrt{u_{k(i)}(1-u_{k(i)})}}\right| \\
+&  \max_{i=1,\ldots,n} \left|\frac{x_i}{\sqrt{x_i(1-x_i)}} - \frac{u_{k(i)}}{\sqrt{u_{k(i)}(1-u_{k(i)})}}\right| \\
\leq &  \qquad 2/n^{9/4} + 4/n^{9/4}. 
\end{align*}
where 
the last step holds for $n > n_0$ by the next lemma.
Reinstating the original $\sqrt{n}$ normalization gives the bound 
\[
 \max_{i=1,\ldots,n}  \sqrt{n}\frac{F(x_i)-x_i}{\sqrt{x_i(1-x_i)}} - \max_{u_k \in \Ucal_n}  \sqrt{n} \frac{F(u_{k})-u_{k}}{\sqrt{u_{k}(1-u_{k})}} \leq  n^{1/2} \cdot \frac{6}{n^{9/4}}, \qquad n > n_0. 
\]
As soon as
$n > \max(n_0,6^{4/3})$, our main claim follows, 
as then $n^{-1} >  6/n^{7/4}$.}
\end{proof}

{\begin{lemma}
\label{lem:2aa}
Let $(x_i)_{i=1}^n$ denote a numerical sequence
taking values in $(0,1)$, obeying
\[
n^{-3/2} < \min_i x_i  \leq \max_i x_i < 1- n^{-3/2}.
\]
Let, 
for each $i$, $u_k(i)$ denote the closest grid point
from our equispaced grid $\Ucal_{n}$ that is as large or larger. Then
\begin{equation} \label{eq:bound1}
    \left |u_{k(i)}-x_{i} \right | \leq n^{-3}; \qquad i=1,\dots,n.
\end{equation}
For $n > n_0$, and $i=1,\dots,n$,
\begin{equation} \label{eq:bound2}
   \left | \frac{1}{\sqrt{u_{k(i)}(1-u_{k(i)})}} - \frac{1}{\sqrt{x_i(1-x_i)}} \right | \leq 2 \cdot n^{-9/4}
\end{equation}
\begin{equation} \label{eq:bound3}
   \left | \sqrt{\frac{u_{k(i)}}{1-u_{k(i)}}} - \sqrt{\frac{x_i}{1-x_i}} \right | \leq 4 \cdot n^{-9/4} 
\end{equation}
\end{lemma}}

\begin{proof}
{\eqref{eq:bound1} follows essentially by definition.}

{The function $w(x) = \frac{1}{\sqrt{x(1-x)}}$ is convex on $(0,1)$,
with derivative $w'(x) = \frac{(1/2-x)}{(x(1-x))^{3/2}}$. 
Now $|w'(x)| \leq (d(x,\{0,1\})/2)^{-3/2}/2$.
Since  $\min_i d(x_i,\{0,1\}) \geq n^{-3/2}$, we get
\[
\max_{i} |w'(x_i)|  \leq (n^{-3/2}/2)^{-3/2}/2 = n^{3/4} \cdot \sqrt{2}.
\]
On $0 < x < u < 1/2$, $w$ is convex decreasing,
so $w(x) > w(u) > w(x) + w'(x)|u-x|$; so there,
 $|w(x)-w(u)| \leq (\max_i |w'(x_i)|) \cdot |u -x|$.
Since $|u-x| \leq n^{-3}$ we have
\[
|w(x)-w(u)| \leq \sqrt{2} \cdot n^{3/4} \cdot n^{-3} = \sqrt{2} \cdot n^{-9/4}; \quad 0 < x < u < 1/2.
\]
Arguing similarly for $1/2 \leq x < u < 1$, where $w$ is convex increasing,
$ w(u) - w'(u)|u-x| \leq w(x) < w(u) $.  To determine
$\max_{1/2 < x_i < 1} |w'(u_{k(i)})|$, use
$|w'(u)| \leq (d(u,\{0,1\})/2)^{-3/2}/2$ and notice that
$\min_i d(u_{k(i)},\{0,1\}) \geq n^{-3/2} - n^{-3} > n^{-3/2}/2$; hence
\[
\max_{1/2 < x_i < 1} |w'(u_{k(i)})| \leq (n^{-3/2}/4)^{-3/2}/2 \leq n^{3/4} \cdot 2,
\]
since $4^{3/2}/2 = 2$. Combining  
$|w(x)-w(u)| \leq (\max_i |w'(u_{k(i)})|) \cdot |u -x|$
and $|u-x| \leq n^{-3}$,
 \eqref{eq:bound2} follows for the case $1/2 \leq x < u < 1$.  
 We omit the argument for the case
$1/2-n^{-1/3} \leq x \leq 1/2 \leq u \leq 1/2+n^{-1/3}$.
\eqref{eq:bound2} follows.}

{Now consider the function $h(x) = \sqrt{\frac{x}{1-x}}$; 
this is increasing  on $(0,1)$,
with derivative $h'(x) = x^{-1/2} (1-x)^{-3/2}/2$. 
Now $|h'(x)| \leq (d(x,\{0,1\})/2)^{-3/2}$.
Since  $\min_i d(x_i,\{0,1\}) \geq n^{-3/2}$, we get
$\max_i |h'(x_i)| \leq n^{-3/4} \cdot 2^{3/2}$.}

{On $0 < x < u < 1/2$, $h$ is concave, so 
$h(x) < h(u) < h(x) + h'(x)|u-x|$; so there,
 $|h(x)-h(u)| \leq (\max_i |h'(x_i)|) \cdot |u -x|$;
since $|u-x| \leq n^{-3}$ we get
\[
|h(x)-h(u)| \leq 2^{3/2} \cdot n^{3/4} \cdot n^{-3} = 4 n^{-9/4}; \qquad 0 < x < u < 1/2.
\]}

{On $1/2 \leq x < u < 1$, $h$ is convex, and 
$ h(u) - h'(u)|u-x| \leq h(x) < h(u) $;
we have the bound
$|h(x)-h(u)| \leq (\max_i |h'(u_{k(i)})|) \cdot |u -x|$,
which requires us to  bound
$\max_{1/2 < x_i < 1} |h'(u_{k(i)})|$.
We notice that
$\min_i d(u_{k(i)},\{0,1\}) \geq n^{-3/2} - n^{-3} > n^{-3/2}/2$; hence
\[
\max_{1/2 < x_i < 1} |h'(u_{k(i)})| \leq (n^{-3/2}/4)^{-3/2}/2 \leq n^{3/4} \cdot 2,
\]
\-- again since $4^{3/2}/2 = 2$. Again using
$|u-x| \leq n^{-3}$, we get
\[
|h(x)-h(u)| \leq 2^{3/2} \cdot n^{3/4} \cdot n^{-3} \leq 4 n^{-9/4}; \qquad 1/2 < x < u < 1-n^{-3}.
\]
 We again omit the argument for the case
$1/2-n^{-1/3} \leq x \leq 1/2 \leq u \leq 1/2+n^{-1/3}$.
\eqref{eq:bound3} follows.}
\end{proof}

\begin{lemma}
\label{lem:Ftilde0}
For any $\eta>0$ and $a>0$, uniformly in $u\in [n^{\eta-1},1/2]$, 
\begin{align*}
    \Pr \left( \sqrt{n} \left( \tilde{F}^{(0)}_n(u) - \ex{\tilde{F}^{(0)}_n(u)} \right) > a \sqrt{u(1-u)} \right) \leq e^{-{a\cdot c} \cdot n^{\eta/2}} + e^{-n^{1-(\beta+1/2)/2}},
\end{align*}
{for some $c=c(\beta) > 0$}.
\end{lemma}
\subsubsection{Proof of Lemma~\ref{lem:Ftilde0}}
\begin{align}
    \sqrt{n}\left(\tilde{F}^{(0)}(u) - \ex{\tilde{F}^{(0)}(u)} \right) = \sqrt{n} \frac{M}{n} \left( H^{(0)}(u)-u\right), 
\end{align}
so, {conditional on $M$,}
\begin{align}
\Pr & \left(\sqrt{n}\left(\tilde{F}^{(0)}(u) - \ex{\tilde{F}^{(0)}(u)} \right) > a \sqrt{u(1-u)}|M \right) \\ 
& {=} \Pr \left( H^{(0)}(u) - u > a \frac{\sqrt{n}}{M} \sqrt{u(1-u)}   | M \right) \\
& = \Pr \left( H^{(0)}(u) > (1+\kappa)u |M \right) =: P_{n,M},
\end{align}
say, where 
\begin{align}
    \kappa := \kappa(u; n,M,a) := a \frac{\sqrt{n}}{M} \sqrt{\frac{1-u}{u}}.
\end{align}
{Conditional on $M$}, 
apply the Chernoff inequality \cite[Ch. 4]{mitzenmacher2017probability}:
\begin{align}
\label{eq:Chernof1}
    P_{M,n} \le \exp\left(- M \frac{\kappa^2}{1+\kappa} u \right).
\end{align}
Consider now the event
$\Omega_n \equiv \{M < n^{1-\beta+\eta_1} \}$. 
{From Lemma~\ref{lem:M}, for any $\eta_1 >0$, 
there exists $n_0 \equiv n_0(\beta,c',\eta_1)$ such that 
\[
\Pr(\Omega_n) \leq c' e^{{-} n^{1-\beta+\eta_1}},\qquad n \geq n_0.
\]
Since $\beta \in (1/2,1)$ we
may specifically choose 
$\eta_1 = (\beta-1/2)/2 >0$, 
yielding $1-\beta+ \eta_1 = 1- (\beta+1/2)/2 > 0$.
%Hence $P(\Omega_n^c) \goto 0$.
}

At the same time, this choice of $\eta_1$ 
gives us $-1/2+\beta-\eta_1 = (\beta-1/2)/2$ ($>0$, as  $\beta > 1/2$).
We are assuming $u \leq 1/2$ so on the event $\Omega_n$:
\begin{align}
    \kappa = a \frac{\sqrt{n}}{M} \sqrt{\frac{1-u}{u}} \ge 
    a \cdot  n^{-1/2+\beta-\eta_1} {= a } \cdot n^{(\beta-1/2)/2} \to \infty. 
\end{align}
Consequently, $\frac{\kappa^2}{1+\kappa} = \kappa \cdot (1 + o(1))$ and, for large $n$, $M \frac{\kappa^2}{2+\kappa} \geq c'' M \cdot \kappa$ for $c'' = c''(\beta)$. Now $M \kappa u = a \sqrt{n} \sqrt{u(1-u)}$, and since $u < 1/2$,
\[
M \kappa u = a \sqrt{n} \sqrt{u(1-u)} > { a \sqrt{n u/2}}.
\]
By hypothesis, $u \geq n^{\eta-1}$, hence
{
\begin{align}
   M \cdot \frac{\kappa^2}{2+\kappa} \cdot  u \geq a \cdot c''' \cdot n^{\eta/2},
\end{align}
for $n > n'''(\beta)$, where $c''' = c'''(\beta)$. Equation \eqref{eq:Chernof1} 
and $\Omega_n$ together imply that for $n > n'''(\beta)$
\begin{align}
     {P_{M,n} \leq \exp\{ {- {a \cdot c'''} \cdot n^{\eta/2}} \}}.
\end{align}
}
Let $P_n$ denote the probability mentioned in the statement
of the lemma; we conclude that
\begin{align}
     {P_n \leq  \Pr(\Omega_n^c) + P_{M,n}} \leq e^{-n^{1-(\beta+1/2)/2}} + e^{-{a \cdot c'''} \cdot n^{\eta/2}} . 
\end{align}
{as advertised.}
\qed
\begin{lemma}
\label{lem:exp_F_diff}
Fix $t_1<1/2$ and let $t_0 = n^{-1}$. For any $\epsilon>0$, there exists $n_0 \equiv n_0(\epsilon)$ such that, if $r < \rho(\beta)$, then
\begin{align*}
\sup_{t \in (t_0,t_1)} \sqrt{n}\left|\ex{\tilde{F}^{(1)}(t)} - \ex{\tilde{F}^{(0)}(t)} \right| w(t) \leq \epsilon,
\end{align*}
for all $n \geq n_0$.
\end{lemma}

\subsubsection{Proof of Lemma~\ref{lem:exp_F_diff}}
We have $\ex{\tilde{F}^{(0)}(t) |M } = tM/n$ for $t \in (0,1)$. 
Consider $t_n \in (t_0,t_1)$, parameterized as $t_n = n^{-q_n}$, 
where $q_n \in (-\log(t_1)/\log(n),1)$. 
From \eqref{eq:H1_I}, we have
\begin{align}
    \ex{\tilde{F}^{(1)}(t_n)|M} = \frac{M}{n} n^{-\alpha(q,r)(1+o_p(1))},
\end{align}
as $n\to \infty$. {Let $\eta>0$. 
On the event $\Omega_n \equiv \{ M < n^{1-\beta+\eta}\}$, 
using also that $q \geq \alpha(q,r)$, we get that $M/n < n^{-\beta+\eta}$;
since $w(t_n) \leq 2 \cdot n^{q_n/2}$ we have (here $q = q_n$):}
\begin{align}
     & \sqrt{n}\left|\ex{\tilde{F}^{(1)}(t_n)\mid M} - \ex{\tilde{F}^{(0)}(t_n)\mid M} \right| w(t_n) \\ 
     & \leq {\sqrt{n} \cdot \frac{M}{n} 
     \left|n^{-\alpha(q,r)(1+o_p(1))} - n^{-q} \right| \cdot n^{q/2} \cdot 2} \\
     & \leq {n^{1/2} \cdot n^{-\beta+\eta} \cdot n^{-\alpha(q,r)(1+o_p(1))} \cdot n^{q/2} \cdot 2} \\
     & = n^{-\beta+\eta {+ \frac{1+q}{2}} -\alpha(q,r)+o_p(1)\alpha(q,r)}.
     \label{eq:proof:exp_F_diff}
\end{align}
As $(\beta,r)$ satisfy \eqref{eq:rho_cond}, we can choose $\eta>0$ such that
\[
\kappa \equiv \beta + \alpha(q,r) {- \frac{1+q}{2}}  - \eta > 0.
\]
Now, $\alpha(q,r)$ is continuous in $q$ by assumption and hence bounded for $q\in [0,1]$,
so $\alpha(q,r)o_p(1) = o_p(1)$. It follows that, on the event $\Omega_n$,
for all $n$ large enough,
\begin{align}
    & \sup_{t \in (t_0,t_1)}  \sqrt{n}
    \left|\ex{\tilde{F}^{(1)}(t)\mid M} - \ex{\tilde{F}^{(0)}(t)\mid M} \right| w(t) \\ 
    & \leq \max_{q\in[0,1]} n^{-\beta+\eta {+ \frac{1+q}{2}}-\alpha(q,r)+o_p(1)} \\
    & = n^{-\kappa+o_p(1)} \to 0.
\end{align}
{Since $F^{(h)}(t) \in [0,1]$,
\[
| \ex  {\ex{\tilde{F}^{(h)}(t)\mid M} | \Omega_n } 
- \ex{\tilde{F}^{(h)}(t)} | 
\leq 2 \cdot \Pr(\Omega_n^c); \qquad h \in \{0,1\}.
\]
Moreover $\sqrt{n} \cdot w(t_n) \leq n$ 
over the range $t_n \in (n^{-1},t_1)$.
Finally, by Lemma~\ref{lem:M}, $ n \cdot \Pr(\Omega_n^c)\to 0$} at an
exponential rate that only depends on $\eta$.
\qed

\begin{lemma}
\label{lem:Ftilde1}
Consider a fixed $u \in (0,1)$. Then for $a>0$ and $\beta \in (1/2,1)$,
\begin{align}
\Pr \left( \sqrt{n} \left( \tilde{F}_n^{(1)}(u) - \ex{\tilde{F}_n^{(1)}(u)} \right) > a \sqrt{u(1-u}  \right) \leq e^{-a \sqrt{n u (1-u)}} + e^{-n^{1/4}}.
\end{align}
\end{lemma}
\subsubsection*{Proof of Lemma~\ref{lem:Ftilde1}}
Define 
\[
\delta = \delta(u;n,M,{a}) \equiv \frac{a \sqrt{n}  \sqrt{u(1-u)} } {M \cdot \ex{H_M^{(1)}(u)} }.
\]
Conditioning on $M$,
\begin{align*}
   & \Pr \left( \sqrt{n} \left( \tilde{F}_n^{(1)}(u) - \ex{\tilde{F}_n^{(1)}(u)} \right) > a \sqrt{u(1-u)} {\mid M }\right) \\
   &  = \Pr \left( H^{(1)}_M(u) - \ex{H^{(1)}_M(u)} >  \frac{a\sqrt{n}}{M} \sqrt{u(1-u)} \right) \\
   & = \Pr \left( H^{(1)}_M(u) > (1+\delta) \ex{H^{(1)}_M(u)} \right) \leq \exp\left\{-M \frac{\delta^2}{1+\delta} \ex{H^{(1)}_M(u)} \right\},
\end{align*}
The last step follows by the Chernoff inequality:
\begin{align}
    \label{eq:Chernoff}
\Pr \left( X \geq (1+\delta)\mu \right) \leq \left( \frac{e^{-\delta}}{(1+\delta)^{1+\delta}} \right)^\mu \leq e^{-\mu \frac{\delta^2}{1+\delta}},
\end{align}
This is valid for any $\delta>1$,
where $X$ is the sum of $M$ independent Bernoulli random variables with $\mu = \ex{X}$ \cite[Ch. 4]{mitzenmacher2017probability}. 

{Set $\eta_1 = (\beta-1/2)/2$.
Since ${\ex{H^{(h)}_M(u)}} \leq 1$,
on the event $\Omega_n \equiv \{ M \leq n^{1-\beta+\eta_1}\}$,  we have}
\begin{align}
    \delta & \ge \frac{a \cdot \sqrt{n}}{n^{1-\beta+\eta_1}} \frac{ \sqrt{(1-u)u}}{\ex{H^{(h)}_M(u)}}  \\
    & \ge a \cdot n^{\beta-1/2-\eta_1} \sqrt{(1-u)u}.
\end{align}
{Since $\beta>1/2$, our choice of $\eta_1$ 
gives  $\beta-1/2-\eta_1 = (\beta-1/2)/2 >0$,
and so $\delta\to \infty$ with increasing $n$.
Hence $\frac{\delta^2}{1+\delta} = \delta (1+o(1))$, and so for large $n$
\[
  M \frac{\delta^2}{1+\delta} \ex{H^{(1)}_M(u)}  \geq  \frac{1}{2} M \cdot \delta \cdot  \ex{H^{(1)}_M(u)}  = \frac{a}{2} \sqrt{n} \sqrt{u(1-u)}.
\]
Therefore, on $\Omega_n$,  }
\begin{align*}
  \Pr \left( \sqrt{n} \left( \tilde{F}_n^{(1)}(u) - \ex{\tilde{F}_n^{(1)}(u)} \right) > a \sqrt{u(1-u)} {\mid M} \right)  \leq  \exp\left\{ - \frac{a}{2} \sqrt{n}\sqrt{u(1-u)} \right\}.
\end{align*}
Finally, Lemma~\ref{lem:M} implies that $\Pr(\Omega_n) < e^{-n^{1-\beta+\eta_1}}$. 
{Our choice $\eta_1 = (\beta-1/2)/2$ gives $1  - \beta+\eta_1 = 3/4 - \beta/2 > 1/4$ since $1 > \beta$,
and so}
\begin{align}
    \Pr \left( \sqrt{n} \left( \tilde{F}_n^{(1)}(u) - \ex{\tilde{F}_n^{(1)}(u)} \right) > a \sqrt{u(1-u)} \right) & \leq e^{-a\sqrt{n}\sqrt{u(1-u)} } + e^{-n^{(3-2\beta)/4}} \\
    & \leq e^{-a\sqrt{n}\sqrt{u(1-u)} } + e^{-n^{1/4}}
\end{align}
 \qed

\newcommand{\Ueta}{\Ucal^\eta_n}
\begin{lemma}
\label{lem:grid1}
{Fix $0 < \eta < 1/2$.
Using the grid $\Ucal_n$ of \eqref{eq:Ucal_def}, 
define the subgrid $\Ueta \equiv \Ucal_n \cap (n^{-1+\eta},1/2)$. 
For $c>0$, and $h \in {0,1}$,}
\begin{align}
     \Pr \left\{  {\max_{u \in \Ueta}} \sqrt{n} \left( \tilde{F}^{{(h)}}(u) - \ex{\tilde{F}^{{(h)}}(u)}  \right)w(u)\geq c  \right\} \to 0. 
\end{align}
\end{lemma}

\subsubsection*{Proof of Lemma~\ref{lem:grid1}}
{Lemmas~\ref{lem:Ftilde0} and \ref{lem:Ftilde1} each imply
exponential inequalities of the form}
\begin{align}
    \max_{u \in (n^{-1+\eta},1/2)} \Pr \left\{ \sqrt{n} \frac{\tilde{F}^{{(h)}}(u) - \ex{\tilde{F}^{{(h)}}(u)}}{\sqrt{u(1-u)}} \geq c \right\} \leq 
    {2 \cdot e^{- (a+bc) \cdot n^{\delta}}},
\end{align}
{where $a,b \geq 0$, $a + bc > 0$, 
$\delta= \delta(\eta)>0$ \--} valid 
for all $n\geq n_0(c)$. 
Therefore, {applying the union bound across the individual grid points}, we get,
{for $h \in \{0,1\}$},
\begin{align}
    &\Pr \left\{  {\max_{u \in \Ueta}} \sqrt{n} \frac{\tilde{F}^{{(h)}}(u) - \ex{\tilde{F}^{{(h)}}(u)}}{\sqrt{u(1-u)}} \geq c  \right\} \\ 
    & \leq \sum_{{u \in \Ueta}} \Pr \left\{ \sqrt{n} \frac{\tilde{F}^{{(h)}}(u) - \ex{\tilde{F}^{{(h)}}(u)}}{\sqrt{u(1-u)}} \geq c \right\} \\
    & \leq {\card( \Ueta) \cdot 2 e^{- (a+b {c}) n^{\delta}}  
    = O(n^4 e^{- (a+b c) n^{\delta}})} \to 0,
\end{align}
since $\card(\Ucal_n)\leq n^4$.
\qed

%% Proof:
\subsection{Proof of Theorem~\ref{thm:converse}}
Set $t_1 = \gamma_0 n < 1/2$ and $t_0 = \bar{Q}_-$.
{
Fix $\delta>0$. Denote by $\Psi_n$ the event that between each consecutive pair $(Q_{(i)}^{(h)},Q_{(i+1)}^{(h)})$,
there exists at least one grid point. From Lemmas~\ref{lem:E1} and \ref{lem:E2}, we have that $\Pr(\Psi_n)>1-\delta/2$ for all $n\geq n_1(\delta)$ for some $n_1(\delta)$. Conditioning on $\Psi_n$, Lemma~\ref{lem:HC_diff}, implies
\begin{align}
\HC_n^{(1)}-\HC_n^{(0)}  \leq \sup_{t \leq t_1} \sqrt{n}  \Delta_n(t) w(t) = \sup_{t \in [ t_0,t_1]} \sqrt{n}  \Delta_n(t) w(t),
\end{align}
while Lemma~\ref{lem:2a} further implies that 
\begin{align}
    \label{eq:proof1}
\sup_{t \in (t_0, t_1)} \sqrt{n}  \Delta_n(t)w(t) \leq \max_{u_k \in \Ucal_n \cap (t_0,t_1)} \sqrt{n} \Delta_n(u_k) w(u_k) + 2 n^{-1}.
\end{align}}
%The probability of the event under which these assumptions hold is tending to one, as implied by Lemmas~\ref{lem:E1} and \ref{lem:E2}. Consequently we may proceed in the proof while assuming that \eqref{eq:proof1} holds. \par
Write
\begin{align*}
\Delta_n(t) & = \left(\tilde{F}^{(1)}(t)-
\ex{\tilde{F}^{(1)}(t)}
\right) + \left(\ex{\tilde{F}^{(1)}(t)}-
\ex{\tilde{F}^{(0)}(t)}
\right)  \\
& \qquad +  \left(\ex{\tilde{F}^{(0)}(t)}-
\tilde{F}^{(0)}(t)
\right).
\end{align*}
Fix $\eta>0$ and $\epsilon>0$. From Lemma~\ref{lem:Qminus},  $n^{\eta-1}\bar{Q}_- > 1$ with probability at least $1-\delta/2$ for all $n\geq n_2(\eta,\delta)$. Therefore, Lemmas~\ref{lem:Ftilde0} and \ref{lem:Ftilde1}, {via Lemma \ref{lem:grid1},}  
{imply that} 
\begin{align*}
    \max_{u_k \in \Ucal_n \cap (t_0,t_1)} \sqrt{n}\left(\tilde{F}^{(0)}(u_k) - \ex{\tilde{F}^{(0)}(u_k)} \right) w(u_k) \leq \epsilon,
\end{align*}
and 
\begin{align*}
\max_{u_k \in \Ucal_n \cap (t_0,t_1)} \sqrt{n}\left(\tilde{F}^{(1)}(u_k) - \ex{\tilde{F}^{(1)}(u_k)} \right) w(u_k) \leq \epsilon,
\end{align*}
for $n \geq n_3(\eta, \epsilon)$ with probability greater than $1-\delta/2$. 
From Lemma~\ref{lem:exp_F_diff} we have that 
\begin{align*}
\sup_{t \in (t_0,t_1)} \sqrt{n}\left(\ex{\tilde{F}^{(1)}(t)} - \ex{\tilde{F}^{(0)}(t)} \right) w(t) \leq \epsilon,
\end{align*}
for all $n\geq n_3(\epsilon)$. Putting it all together, 
\begin{align}
    \HC_n^{(1)}-\HC_n^{(0)}  \leq 3 \epsilon + 2n^{-1},
\end{align}
with probability at least $1-\delta$
for all $n \geq n_4(\epsilon,\eta,\delta)$. As $\delta$, $\eta$ and $\epsilon$ are arbitrary, the proof is completed. 

\begin{funding}
This work is supported in part by funding from the Koret Foundation and the NSF under Grant No.~DMS-1816114.
\end{funding}

\bibliographystyle{spbasic}      % basic style, author-year citations
\bibliography{HigherCriticism}   % name your BibTeX data base

\end{document}